\documentclass{article}
\usepackage[a4paper, total={6in, 8in}]{geometry}
\usepackage{authblk}

\usepackage{type1cm}        
\usepackage{makeidx}         
\usepackage{graphicx}        
\usepackage{multicol}        
\usepackage[bottom]{footmisc}

\usepackage{amsmath,amssymb,mathtools,amsthm,stmaryrd, latexsym}
\usepackage[linesnumbered,ruled,vlined, onelanguage]{algorithm2e}
\usepackage[noend]{algorithmic}
\usepackage{enumitem}
\usepackage{diagbox}
\usepackage{nicefrac}
\usepackage{url}
\usepackage{cleveref}

\makeatletter
\newcommand{\refcheckize}[1]{%
  \expandafter\let\csname @@\string#1\endcsname#1%
  \expandafter\DeclareRobustCommand\csname relax\string#1\endcsname[1]{%
    \csname @@\string#1\endcsname{##1}\wrtusdrf{##1}}%
  \expandafter\let\expandafter#1\csname relax\string#1\endcsname
}
\makeatother

\theoremstyle{definition}
\newtheorem{asm}{Assumption}[section]

\theoremstyle{remark}
\newtheorem{rk}{Remark}[section]

\newcommand{\Th}{\mathcal{T}_{h}}

\newcommand{\mbf}[1]{\boldsymbol{#1}}

\newcommand{\bb}[1]{\mathbb{#1}}
\newcommand{\vertiii}[1]{{\left\vert\kern-0.25ex\left\vert\kern-0.25ex\left\vert #1
            \right\vert\kern-0.25ex\right\vert\kern-0.25ex\right\vert}}

\newcommand{\twofigswidth}{0.45\textwidth}
\newcommand{\threefigswidth}{0.31\textwidth}
\SetKwInput{KwInit}{Initialisation}  

\usepackage{newtxtext}       %
\usepackage[varvw]{newtxmath}       

\makeindex             

\title{A multigrid and neural network approach to reduce the computational cost of $\varphi$-FEM}
\author[1]{Raphaël Bulle}
\author[1]{Michel Duprez}
\author[2]{Vanessa Lleras}
\author[2]{Killian Vuillemot}
\affil[1]{Université de Strasbourg, CNRS, Inria, ICube, F-67000 Strasbourg, France}
\affil[2]{IMAG, Univ Montpellier, CNRS UMR 5149, 499-554 Rue du Truel, 34090 Montpellier, France}

\date{}
\begin{document}
%
\maketitle

\abstract{In this work, we present a combination of a multigrid approach and the $\varphi$-FEM immersed boundary finite element method to reduce its computational cost while preserving its accuracy. 
To further reduce the numerical cost of the approach, 
we also propose the combination of the previous technique with some neural network methods. 
We illustrate the efficiency of these two approaches with numerical test cases in 2D and 3D.}

\section{Introduction}

The finite element method (FEM) is among the most used methods to solve partial differential equations (PDEs).
However, FEM faces important limitations when considering complex geometries (e.g. modelling of organs in the context of biomedical simulations) since the use of conforming meshes is mandatory to preserve accuracy. To bypass this difficulty, several nonconformal methods have been introduced, such as CutFEM \cite{burman1,burman2}, Shifted Boundary Method \cite{sbm}, or more recently $\varphi$-FEM \cite{phifem,phiFEM2,phifem_chapter, phiFEM_heat}. 
Another class of methods for solving such equations is the Finite Difference method, which relies on structured grids and is thus limited to simple geometries. 
However, in \cite{shortley1938numerical}, a method was proposed to extend the finite difference method to complex geometries. 

To accelerate the computation, the Cartesian structure of the mesh used in the finite difference scheme enables an efficient combination with multigrid approaches \cite{bramble}.
Recently, \cite{phiFD} adapted the $\varphi$-FEM approach to develop a finite difference scheme called $\varphi$-FD and combined it with a multigrid method to solve the Poisson equation with homogeneous Dirichlet boundary conditions. 

Several recent approaches investigate the integration of neural networks within multigrid framework. These methods leverage data-driven models to enhance fine-scale corrections or replace selected components of classical multigrid techniques, while retaining coarse-grid efficiency and theoretical insights inherited from traditional multilevel methods. Hybrid finite element–neural network solvers \cite{richter} and neural multigrid corrections exemplify this trend, aiming to improve accuracy and computational performance for elliptic problems. These hybrid approaches use PINNs \cite{PINN-MG} or deep neural networks \cite{DNN}, but there is not yet any coupling with Fourier Neural Operator (FNO).


To bypass the issue of the computational cost, a combination of $\varphi$-FEM with FNO \cite{paper_FNO}, namely $\varphi$-FEM-FNO \cite{phifem_fno}, has been proposed. 
However, on high complexity problems, such as non-linear elastic problems, the gain in computing time comes at the expense of accuracy since the neural network does not succeed in reaching the precision of FEM-based methods. In this paper, we propose a new approach that has the advantage of preserving the accuracy of the $\varphi$-FEM approach combining neural networks and multigrid methods, while reducing its computation time. 

Consider $\Omega \subset \bb{R}^d$ (here $d=2, 3$) a domain of boundary $\Gamma$ described by a level-set function $\varphi$, i.e. 
     $ \Omega:= \{ \varphi < 0 \} $ and $ \Gamma:= \{ \varphi = 0 \}\,.$
To introduce our method, we consider two equations: on the one hand, the Poisson equation with Dirichlet boundary conditions, 
\begin{equation}
     \label{eq:main_poisson_equation_dirichlet_homo}
            - \Delta u  = f\, \ \hfill \text{ in } \Omega\, \qquad \text{ and } \qquad
            \hfill u    = 0\, \  \hfill \text{ on } \Gamma\,,
\end{equation}
where $f \in L^2(\Omega)$.
On the other hand, a non-linear Poisson equation, 
\begin{equation}
      \label{eq:non_linear_poisson_main}
            - \nabla \cdot (q(u) \nabla u)  = f\, \hfill \text{ in } \Omega\,, \qquad \text{ and } \qquad
            \hfill  u                       = 0\, \hfill \text{ on } \Gamma\,,  
\end{equation}
with $f \in L^2(\Omega)$ and $q(u)$ a non-linear function.

The manuscript is organized as follows:
Section 2 introduces the $\varphi$-FEM scheme. Section 3 is devoted to the combination of this scheme with a multigrid approach and presents the corresponding numerical results. Section 4 proposes a hybrid extension based on neural networks, combining reduced computational cost with guaranteed numerical convergence. Finally, Section 5 concludes the paper.

\section{The $\varphi$-FEM method}

Let us introduce the mathematical setting of the method. We consider the original $\varphi$-FEM scheme \cite{phifem} to solve \cref{eq:main_poisson_equation_dirichlet_homo}. However, this can be extended to all the other $\varphi$-FEM schemes that have been proposed in the literature (see e.g. \cite{phifem_chapter}). The method relies on the use of a level-set function to describe the considered geometry and its boundary, allowing us to rewrite the unknown as $u = \varphi w$ with $w$ the new unknown of the problem. 

We consider $\mathcal{T}_h^\mathcal{O}$ a Cartesian grid covering the domain $\mathcal{O}$, with triangular cells of size $h$. Let $\varphi_h = I_{h, \mathcal{O}}^{(l)} \varphi$ the standard Lagrange continuous interpolation of $\varphi$ (of degree $l>0$) on $\mathcal{T}_h^\mathcal{O}$.
Using this interpolation, we construct the submeshes $\mathcal{T}_h$ and $\mathcal T_h^{\Gamma}$ of $\mathcal{T}_h^\mathcal{O}$ collecting all the cells of the grid intersecting $\{ \varphi_h < 0 \}$ and intersecting $\{\varphi_h = 0\}$ respectively, as well as the subsets covered by these cells.
In other words,
\begin{align}
      \mathcal{T}_h &:= \left\{T \in \mathcal{T}_h^{\mathcal{O}} : \ T \cap \{ \varphi_h < 0 \} \neq \emptyset  \right\},\quad \Omega_h := \bigcup_{T \in \mathcal T_h} \overline T,\label{eq:Th}\\
      \mathcal T_h^{\Gamma} &:= \left\{T \in \mathcal{T}_h^{\mathcal{O}} : \ T \cap \{ \varphi_h = 0 \} \neq \emptyset  \right\},\quad \Omega_h^{\Gamma} := \bigcup_{T \in \mathcal T_h^{\Gamma}} \overline T.
      \label{eq:Th_Gamma}
\end{align}
In addition, we construct a set containing facets of $\mathcal T_h^{\Gamma}$,
\begin{equation}
    \mathcal{F}_h^\Gamma:= \{ F \in \mathcal{T}_h^\Gamma \setminus \partial\Omega_h  \} \,.
    \label{eq:F_h_Gamma}
\end{equation}
Let us finally introduce the considered finite element spaces: for $k\geqslant1$,
\begin{equation*}
      V_h^{(k)}:= \{ v_h \in H^1(\Omega_h) : \ v_h|_T \in \bb{P}_k(T) \ \forall \ T \in \mathcal{T}_h \} \,.
\end{equation*}
The $\varphi$-FEM scheme to solve \cref{eq:main_poisson_equation_dirichlet_homo} is then given by  find $w_h \in V_h^{(k)}$ for all $s_h \in V_h^{(k)}$, such that denoting $u_h = \varphi_h w_h$, and $v_h = \varphi_h s_h$,
\begin{equation}\label{eq:scheme_poisson_dirichlet_phi_fem}
    \int_{\Omega_h} \nabla u_h \cdot \nabla v_h - \int_{\partial\Omega_h} \frac{\partial u_h}{\partial n} v_h
    + G_h^{lhs}(u_h,v_h)
    = \int_{\Omega_h} f v_h + G_h^{rhs}(v_h)\,,
\end{equation}
with
\begin{equation*}
      G_h^{lhs}(u,v) = \sigma_D h \sum_{E \in \mathcal{F}_h^\Gamma} \int_E \left[ \frac{\partial u}{\partial n} \right]\left[ \frac{\partial v}{\partial n} \right]  + \sigma_D h^2 \sum_{T \in \mathcal{T}_h^\Gamma}\int_T \Delta u \Delta v \,,
\end{equation*}
and 
 $ \qquad \qquad  \qquad \qquad   G_h^{rhs}(v) = - \sigma_D h^2 \sum_{T \in \mathcal{T}_h^\Gamma}\int_T f \Delta v \,.$

\begin{rk}\label{rk:non_homo_dirichlet}
In the case of non-homogeneous Dirichlet boundary conditions $u = g$ on $\Gamma$, the unknown $u_h$ is written as $u_h= \varphi_h w_h + g_h$, with $g_h$ an interpolation of $g$ on $V_h^{(k)}$ and the scheme is adapted according to this change. 
\end{rk}

\section{The method $\varphi$-FEM-Multigrid}

\subsection{Methodology}

Let us now describe more precisely our method in the case of \cref{eq:scheme_poisson_dirichlet_phi_fem}.
In the following, we denote $V_h^{(k,i)}$ the Lagrange finite element space associated to the mesh $\mathcal T_h^{(i)}$.
The ``coarse step'' starts with an initial cartesian grid $\mathcal{T}_h^{\mathcal{O}, (0)}$. 
We compute the interpolation $\varphi_h^{(0)} \in V_h^{(k,0)}$ of $\varphi$ and use this interpolation to construct $\Omega_h^{(0)}$, $\Th^{(0)}$, $\Omega_h^{\Gamma,(0)}$ and $\Th^{\Gamma,(0)}$ as described in \cref{eq:Th,eq:Th_Gamma} and $\mathcal F_h^{\Gamma, (0)}$ following \cref{eq:F_h_Gamma}.
Thus, we compute $u_{h,(0)} \in V_h^{(k,0)}$, solution to \cref{eq:scheme_poisson_dirichlet_phi_fem} on $\Omega_h^{(0)}$.

Then, we iteratively perform the ``fine steps''.
We construct a sequence of meshes $\{\Th^{(i)},\ i=1,\ldots\}$, where for each $i$, $\Th^{(i)}$ is given by \cref{eq:Th} replacing $\Th^O$ by $\Th^{\mathcal O, (i)}$, a uniform refinement of $\Th^{(i-1)}$ and $\varphi_h$ by $\varphi_h^{(i)} \in V_h^{(k,i)}$ the Lagrange interpolation of $\varphi$.
Thus, $\Th^{\Gamma,(i)}$, $\Omega_h^{(i)}$ and $\Omega_h^{\Gamma, (i)}$ are constructed in a similar way to \cref{eq:Th,eq:Th_Gamma,eq:F_h_Gamma}.
Then, we interpolate $u_{h,(i-1)}$ onto the finite element space $V_h^{(k,i)}$ and use this interpolation to initialize the iterative solver used to compute $u_{h,(i)}$, solution to \cref{eq:scheme_poisson_dirichlet_phi_fem} in $V_h^{(k,i)}$.

We would like to emphasize that, depending on the fineness of the initial mesh $\Th^{(0)}$, some parts of $\Omega$ might be lost and never recovered during the iterative process since for each $i=1,\ldots$, we have $\Omega_h^{(i)} \subset \Omega_h^{(i-1)}$.
A way to avoid this is using the following assumption:
\begin{asm}
    We assume that the exact domain $\Omega$ is a subset of $\Omega_h^{(0)}$.
    \label{asm:exact_domain}
\end{asm}

\begin{figure}
      \centering
      \includegraphics[width=\textwidth]{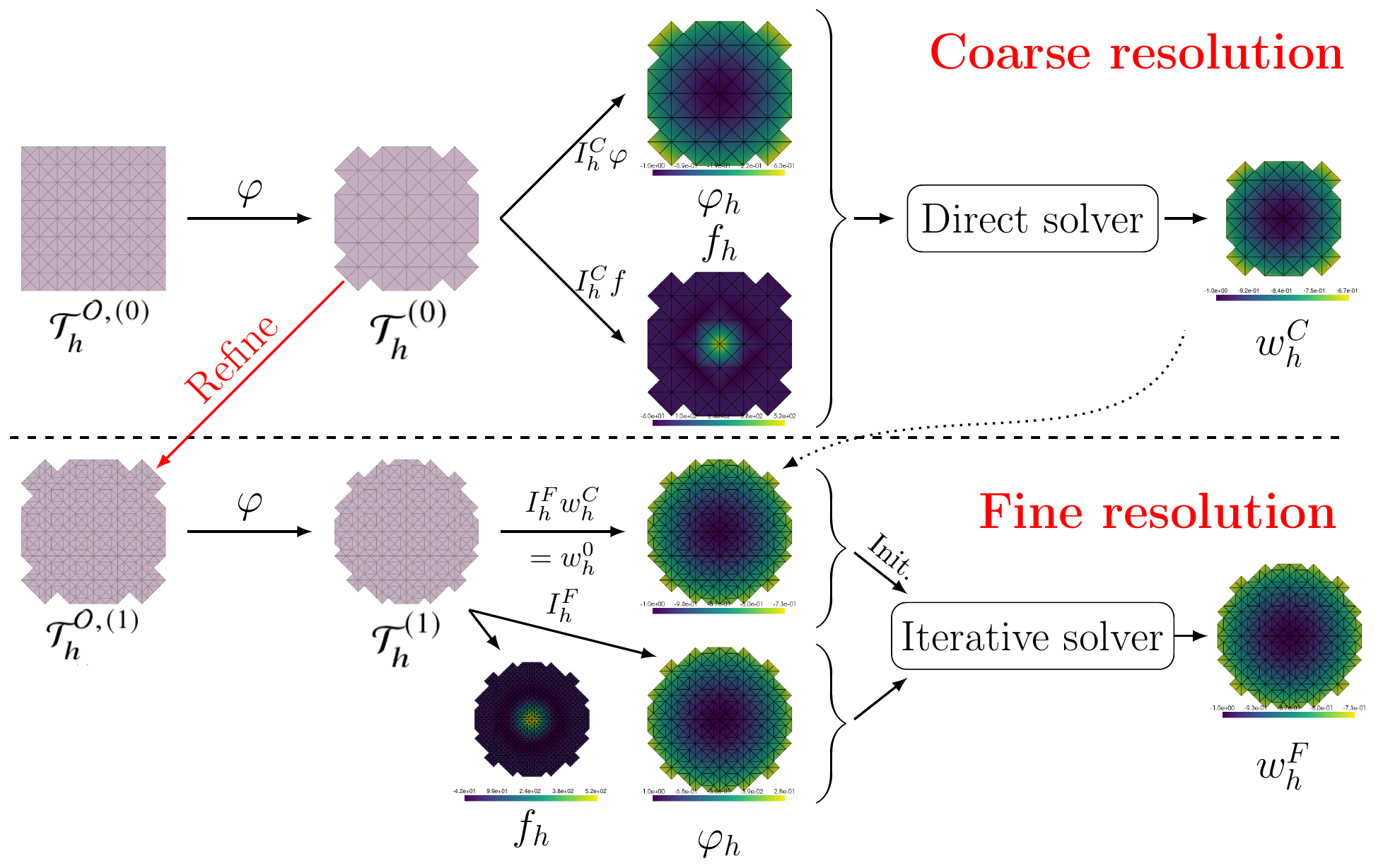}
      \caption{Graphical representation of the pipeline of the $\varphi$-FEM-Multigrid approach. 
      }\label{fig:pipeline_phi_femm}
\end{figure}

\Cref{alg:phifem_multigrid_v2} summarizes our approach in the case of a single ``fine step'' but can be straightforwardly extended to several ``fine steps''.
A graphical representation of the pipeline considered in \cref{alg:phifem_multigrid_v2} is given in \cref{fig:pipeline_phi_femm}, for the case of \cref{eq:main_poisson_equation_dirichlet_homo}.
The Lagrange interpolation operators are denoted by $I_h^{\bullet}$ with $\bullet=C,F$ when the interpolation is performed on the coarse and fine meshes respectively. Similarly, $w_h^{\bullet}$ represents the $\varphi$-FEM solution on the coarse mesh when $\bullet=C$ and on the fine mesh when $\bullet=F$.
\begin{figure}
      \centering
      \includegraphics[width=\threefigswidth]{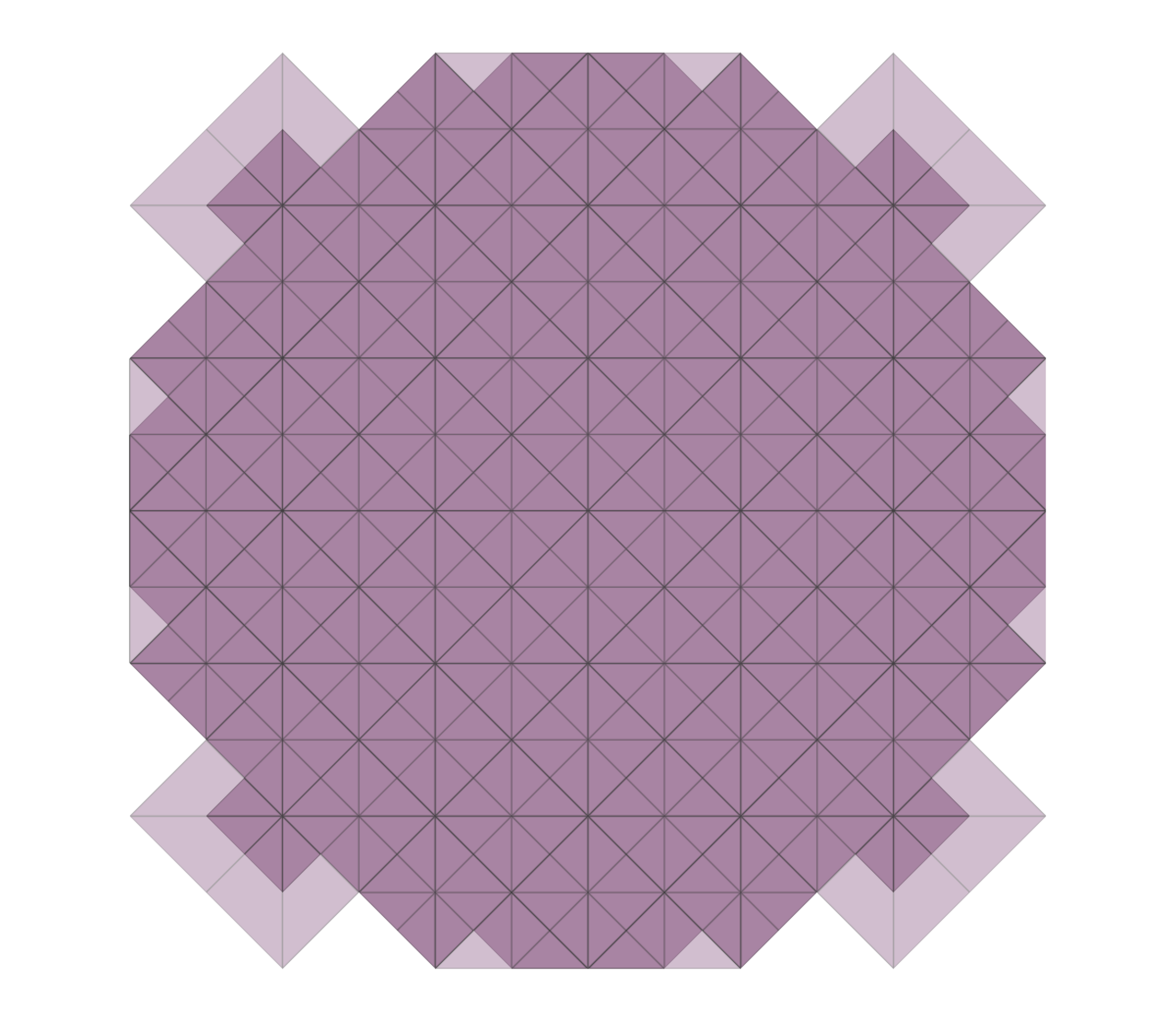}
      \caption{Construction of the $\Th^{(i)}$ mesh (in pink) from the $\Th^{(i-1)}$ mesh (shaded).}\label{fig:phi_fem_multigrid_raffinement_v2}
\end{figure}

                        


                       
      
      \begin{algorithm}
      \caption{{$\varphi$}-FEM Multigrid}\label{alg:phifem_multigrid_v2}
        \KwIn{$N$ : number of refinement steps, $n$ : number of cells in each direction of the coarse grid.}
        
        Generate $\Th^{\mathcal{O},(0)}$, with $n \times n$ cells

        Interpolate $\varphi$ into $V_h^{(k,0)}$ to get $\varphi_h^{(0)}$

        Compute $\Th^{(0)}$, $\Th^{\Gamma,(0)}$, $\mathcal{F}_h^{\Gamma, (0)}$

        If using an iterative solver, initialize with $0$
        
        Solve \cref{eq:scheme_poisson_dirichlet_phi_fem} on $\Th^{(0)}$ to get $u_{h,(0)} \in V_h^{(k,0)}$
  
        $\Th^{\mathcal{O},(1)} \leftarrow \Th^{(0)}$

        Interpolate $\varphi$ into $V_h^{(k,1)}$ to get $\varphi_h^{(1)}$

        Compute $\Th^{(1)}$, $\Th^{\Gamma,(1)}$, $\mathcal{F}_h^{\Gamma, (1)}$
        
        Interpolate $u_{h,(0)}$ into $V_h^{(k,1)}$
        
        Initialize the iterative solver with $u_{h,(0)}$

        Solve \cref{eq:scheme_poisson_dirichlet_phi_fem} on $\Th^{(1)}$ to get $u_{h,(1)} \in V_h^{(k,1)}$.
      \end{algorithm}

\begin{rk}
The mesh refinement is performed by the \texttt{refine} routine of dolfinx \cite{dolfinx}, which consists in adding a new node at the midpoints of each edge of the mesh.
\end{rk}

\subsection{Numerical results}

Let us now illustrate the interest of the approach on two numerical results, by comparing the method to the standard conforming finite element method (referred to as "Standard-FEM") and to the $\varphi$-FEM technique. The numerical simulations are performed in Python using dolfinx \cite{dolfinx} (version 0.8.0), and the implementation is available on the following GitHub repository: 
    \url{https://github.com/PhiFEM/PhiFEM_Multigrid.git}.
The simulations were executed on a laptop with an \texttt{Intel Core i7-12700H CPU}, $32$Gb of memory, and an \texttt{NVIDIA RTX A2000 GPU} with $8$Gb of memory.


\subsubsection{Test case 1: a non-linear equation on a disk}

We first consider the case of \cref{eq:non_linear_poisson_main}, where $\Omega$ is a disk centered in $(0.5, 0.5)$ of radius $\nicefrac{\sqrt{2}}{4}$. We take $q(u)=1 + u^3 \exp(2.5 u)$ and $f$ so that the exact solution reads $ u(x,y) = \cos \left( \frac{\pi}{2} r(x,y) \right) $
with $r(x,y) = \frac{1}{R}\sqrt{{(x-0.5)^2} + {(y-0.5)^2}}$. The non-linear function is given by $q(u) = 1 + u^3 \exp(2.5 u)$ and $f$ is computed analytically. To solve the problem, we adapt the $\varphi$-FEM scheme \cref{eq:scheme_poisson_dirichlet_phi_fem} to the case of the \cref{eq:non_linear_poisson_main} by multiplying the relevant terms by $q(u)$.



\begin{figure}
      \centering
      \includegraphics[width=0.8\textwidth]{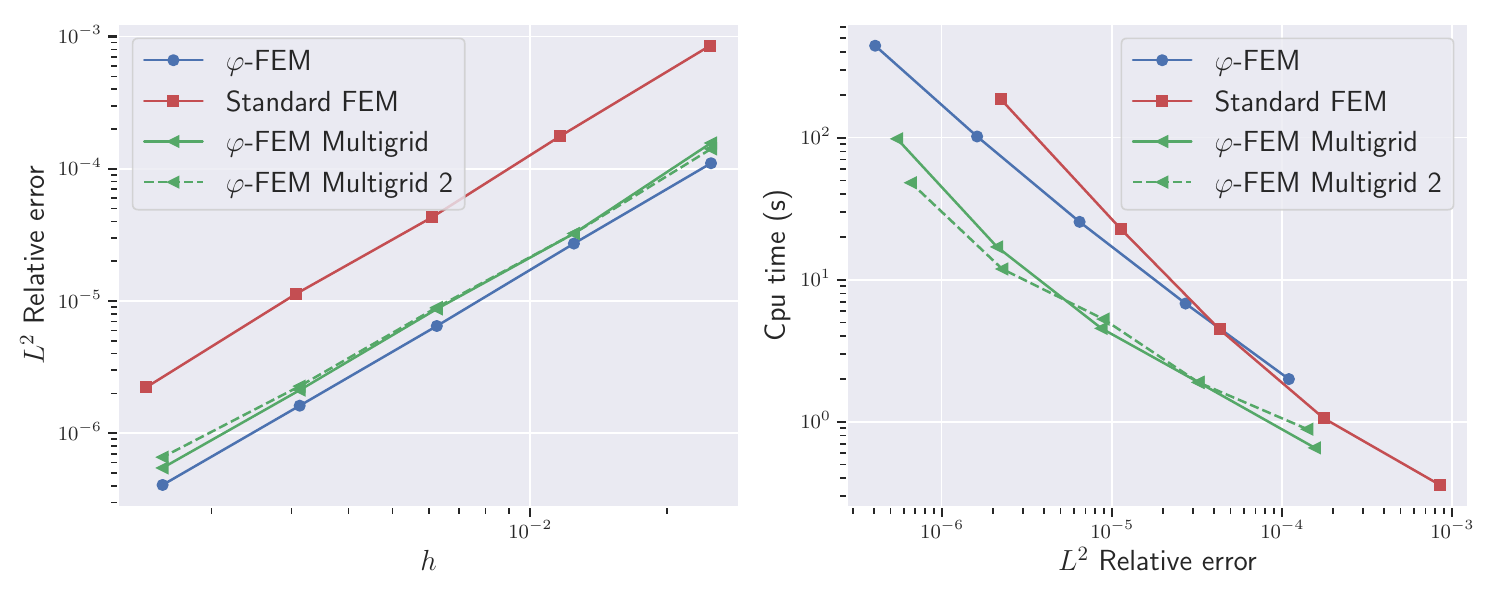}\caption{\textbf{Test case 1.} Left: $L^2$ relative errors with respect to the cell size. Right: CPU time (in seconds) with respect to the $L^2$ relative errors.
      }\label{fig:phi_femm_results_poisson_non_linear}
\end{figure}

The results are given in \cref{fig:phi_femm_results_poisson_non_linear}. In this test case, we propose two versions of the multigrid approach: 
\begin{itemize}
    \item for the first one, coined ``$\varphi$-FEM-Multigrid'', a single ``fine step'' is performed (therefore the resolution of the meshes on the ``coarse step'' changes), and we apply two refinement steps to the coarse grid to obtain the fine grid (without any intermediary resolution step);
    \item for the second version, referred to as ``$\varphi$-FEM-Multigrid 2'', we keep the coarse resolution constant (a mesh with $20 \times 20$ cells) and increase the number of refinement steps. This implementation is the foundation of the combination with neural networks that is proposed later on.
\end{itemize}
Moreover, the tolerance of the iterative solvers is set to $10^{-5}$. In \cref{fig:phi_femm_results_poisson_non_linear}, we observe that both $\varphi$-FEM Multigrid approaches are close to reaching the precision of the classical $\varphi$-FEM method, much faster than both $\varphi$-FEM and Standard-FEM. 

\subsubsection{Test case 2: Poisson equation in a sphere}

\begin{figure}
      \centering
      \includegraphics[width=0.8\textwidth]{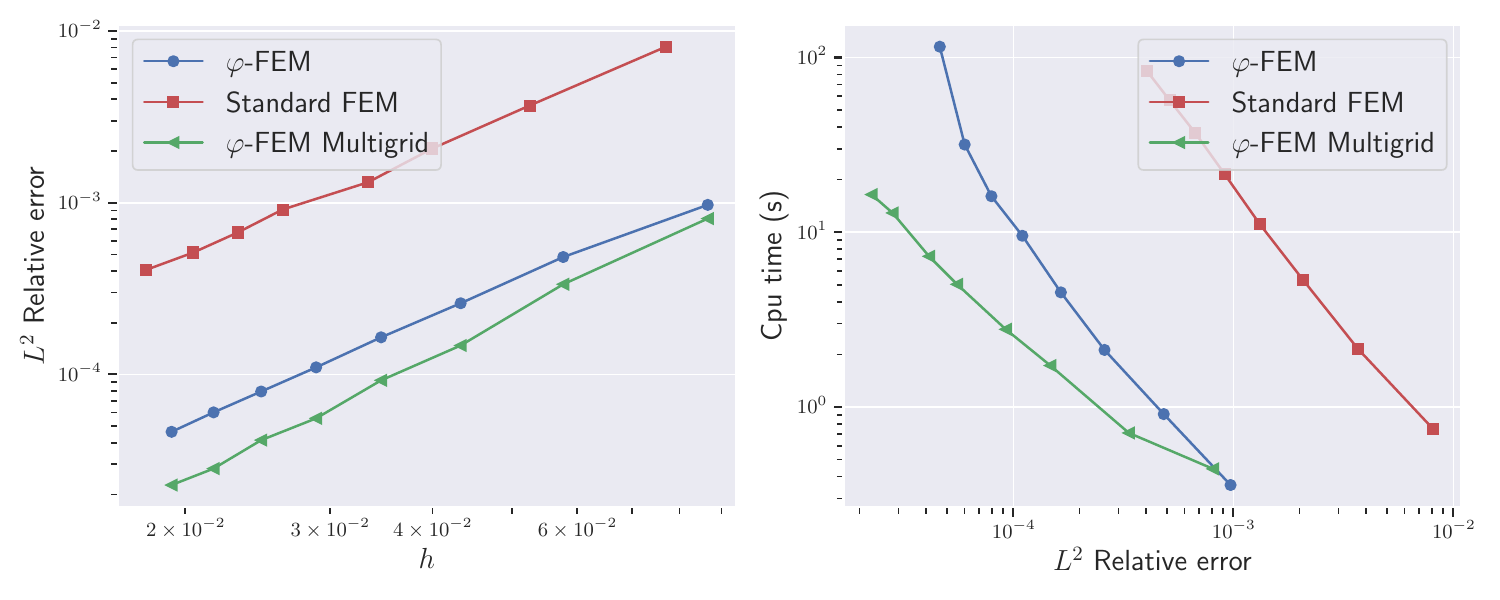}\caption{\textbf{Test case 2.} Left: $L^2$ relative errors with respect to the cell size $h$. Right: CPU time (in seconds) with respect to the $L^2$ relative errors.}\label{fig:phi_femm_results_poisson_dirichlet_3D}
\end{figure}

Let us now extend our study to a 3D geometry of sphere centered in $(0.5, 0.5, 0.5)$, with radius $R = \sqrt{2}/4$. We now consider the \cref{eq:main_poisson_equation_dirichlet_homo}, with the same radial solution, extended to the 3D case, hence with $r =  \frac{1}{R}\sqrt{{(x-0.5)^2} + {(y-0.5)^2}  + {(z-0.5)^2}}$. In this 3D test case, solving the scheme  \cref{eq:scheme_poisson_dirichlet_phi_fem} using a direct solver is prohibitively expansive. Hence, for the classical $\varphi$-FEM approach as for the Standard FEM, we use an iterative solver, namely the BiConjugate Gradient Stabilized. The results of the three methods on such a test case are given in \cref{fig:phi_femm_results_poisson_dirichlet_3D}, highlighting the interest of our approach since we achieve better precision faster than the two other considered methods. 
Note that we do not compare the methods with $\varphi$-FEM-Multigrid-2 in this test case since this would involve using only 3 points on the convergence curves due to the quick increase of computational cost.

\section{{$\varphi$}-FEM-M-FNO: combining multigrid and neural networks}

Let us now introduce a natural extension to the previously proposed approach. Indeed, in \cite{phifem_fno} we introduced $\varphi$-FEM-FNO, a combination between $\varphi$-FEM and the Fourier Neural Operator \cite{paper_FNO}, to bypass the limitation of the FNO to simple geometries. 

\subsection{Methodology}

The $\varphi$-FEM-FNO technique relies on encoding the geometry using the level-set function $\varphi$ and uses this function as an input to a neural operator along with the functions defining the problem to solve. For example, considering the Poisson-Dirichlet \cref{eq:main_poisson_equation_dirichlet_homo}, the operator to train would be given by, 
\begin{align*}
    \mathcal{G}_\theta : \quad \bb{R}^{n_x \times n_y \times 3} \quad &\to \quad \bb{R}^{n_x \times n_y \times 1}\,, \\ 
                                (f_h, g_h, \varphi_h) \quad &\mapsto \quad w_\theta\,,
    \end{align*}
with $w_\theta$ an approximation of the $w_h$ solution obtained using $\varphi$-FEM. 

Then, combining this approach with the use of the multigrid technique seems natural. We propose the following methodology: we replace the first coarse finite element resolution by a call to a pre-trained $\varphi$-FEM-FNO on coarse data.
This multigrid approach, not in real-time (contrary to the classical $\varphi$-FEM-FNO), has several advantages:
\begin{itemize}
    \item the data generation could be lighter than for the application of the classical $\varphi$-FEM-FNO since the multigrid approach requires less data, with coarser resolution;
    \item the training cost is also reduced, due to the lighter dataset. Moreover, since the operator does not need to match high accuracy, the number of training steps can also be reduced, also decreasing the training cost.
    \item The numerical cost is lighter than that of $\varphi$-FEM-Multigrid since the first coarse resolution is avoided, as all the interpolations on the coarse mesh. 
\end{itemize}

However, it is important to highlight that a mandatory, costly step arises with this method. Indeed, a predicted solution obtained with $\varphi$-FEM-FNO is simply a matrix (or a tensor) that would need to be discretized into a finite element function. This transformation is done using the matrix values, which by definition are the nodal values of a finite element function, and constructing a mapping between the matrix indices and the nodal indices. This step is at this time not optimal under the dolfinx \cite{dolfinx} framework, but can be constructed only once for a fixed resolution. 

\begin{rk}
It would be possible to predict the fine solution of the problem directly and use it as an initial solution of the iterative solver. This would avoid the use of the multigrid technique. We find this approach, for example, in \cite{odot:hal-03327818}, where the neural network prediction is used to initialize a Newton solver. However, the data generation (if using a data-based method) and the training steps being done at the fine resolution induce a higher computational cost. In our case, this would also induce the need to construct directly the fine $\Th$ mesh, which is much more costly.

Another approach combining finite element methods and neural networks \cite{correction_Fred} showed very interesting theoretical and numerical results. This method uses a pre-trained PINN of a large number of collocation points, which is corrected by a coarse conforming finite element method. 
\end{rk}

\bigskip

\begin{figure}[!ht]
      \centering
      \includegraphics[width=\textwidth]{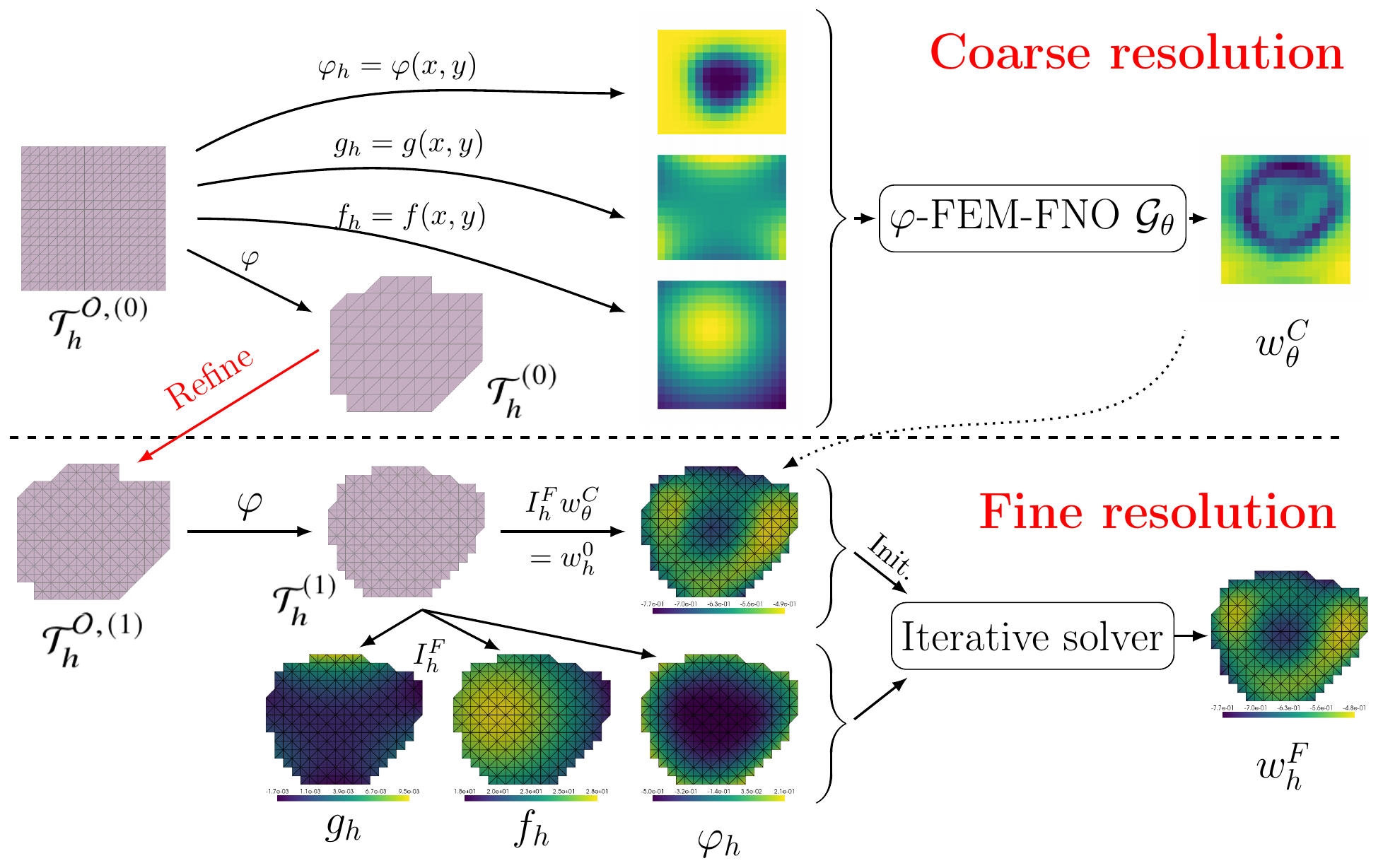}
      \caption{Pipeline of $\varphi$-FEM-M-FNO to solve \cref{eq:non_linear_poisson_main} with non-homogeneous boundary conditions. 
      }\label{fig:pipeline_phi_femm_fno}
\end{figure}

A graphical representation of the pipeline is given in \cref{fig:pipeline_phi_femm_fno}, in the case of the resolution of the \cref{eq:non_linear_poisson_main} with non-homogeneous Dirichlet boundary conditions (i.e. $u=g$ on $\Gamma$). \\ 
The approach consists of three main steps:
\begin{itemize}[label=\textbullet]
      \item Coarse resolution: prediction of a coarse approximation thanks to the \cref{rk:non_homo_dirichlet} \textbf{and} construction of the mesh $\Th^{(0)}$;
      \item Refinement: refinement loop to reach the desired resolution;
      \item  Fine resolution: interpolation of the coarse approximation on the fine mesh and $\varphi$-FEM
      resolution with an iterative solver with a non-zero initialization. 
\end{itemize}













\subsection{Numerical test cases}

\subsubsection{Test case 3: the 2D case}

\begin{figure}[!ht]
      \centering
      \includegraphics[width=\textwidth]{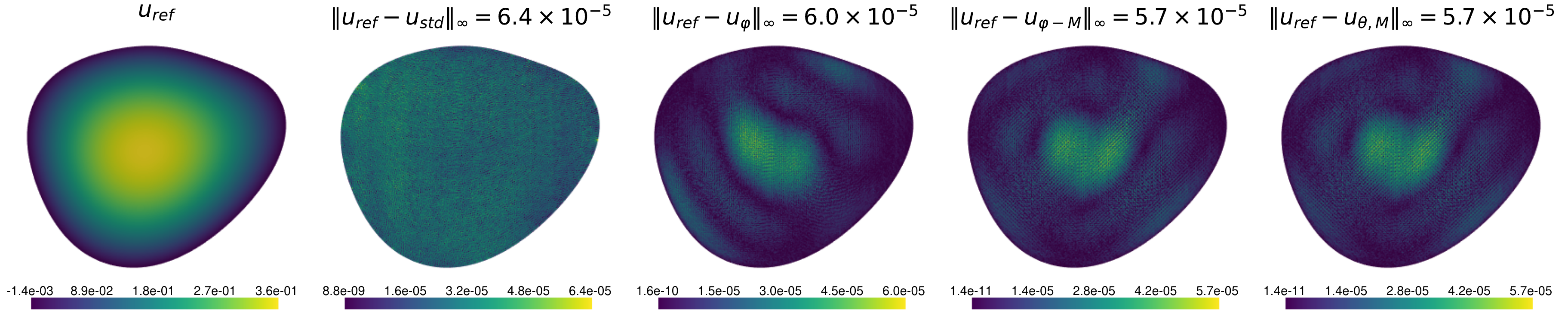}
      \caption{\textbf{Test case 3.} From left to right: reference solution, then difference between the reference solution and the projection of the Standard-FEM solution ($u_{\text{std}}$), the  $\varphi$-FEM solution ($u_{\varphi}$), the $\varphi$-FEM-Multigrid solution ($u_{\varphi-M}$), and the $\varphi$-FEM-M-FNO solution ($u_{\theta,M}$).}
      \label{fig:outputs_phi_fem_fno_multigrid_cas_test_1}
\end{figure}

In this test case, we consider the \cref{eq:non_linear_poisson_main} with non-homogeneous Dirichlet boundary conditions, $u=g$ on $\Gamma$, with 
   $ g_{(\alpha, \beta)}(x,y) = \alpha \left((x-0.5)^2 - (y-0.5)^2\right) \cos\left( \beta y \pi \right)\,,$
where $\alpha$, $\beta \sim \mathcal{U}([-0.8,0.8])$. The considered geometries are defined by level-set functions 
\begin{equation*}
      \varphi(x,y) = (-1)^5 \prod_{j=1}^5 \bigg( -1 + \exp\left(- \frac{x_j^2}{2l_{x,j}^2} - \frac{y_j^2}{2l_{y,j}^2} \right)  \bigg) - 0.5\,,
\end{equation*}
with, for some parameters $(x_{0,j}, y_{0,j}, l_{x,j}, l_{y,j}, \theta_j)$ randomly chosen, $ x_j = \cos(\theta_j) (x-x_{0,j}) - \sin(\theta_j) (y-y_{0,j})$ and $ y_j = \sin(\theta_j) (x-x_{0,j}) + \cos(\theta_j) (y-y_{0, j})$
. The bounds of the random parameters of $\varphi$ are fixed to ensure that the geometries lie inside $[0,1]^2$ and are connected. The FNO is trained with 500 training data and 300 validation data, for 2000 training epochs. The loss function to be minimized is an approximation of the $H^1$ error to the $\varphi$-FEM solution on the considered geometry (cf. \cite[Paragraph 3 and equation (4)]{phifem_fno} for more details).
The method is compared to the Standard-FEM, $\varphi$-FEM, and $\varphi$-FEM-Multigrid. The reference solutions used to compute the errors are obtained by fine Standard-FEM resolutions. An example of output for new data (unseen during training of the FNO) is given in \cref{fig:outputs_phi_fem_fno_multigrid_cas_test_1}. 

To illustrate the performance of the method, we consider two situations, corresponding to two different coarse resolutions. For both, we consider 5 data taken from a test dataset.

\paragraph*{$\mathbf{16 \times 16}$ grids}
We first consider data generated on $16 \times 16$ grids. For a fair comparison, the coarse resolutions of $\varphi$-FEM-Multigrid are performed with the same grid size. We then compare the results at some fine sizes:  32, 64, 128, and 256. 

\begin{figure}[!ht]
      \centering
      \includegraphics[width=0.8\textwidth]{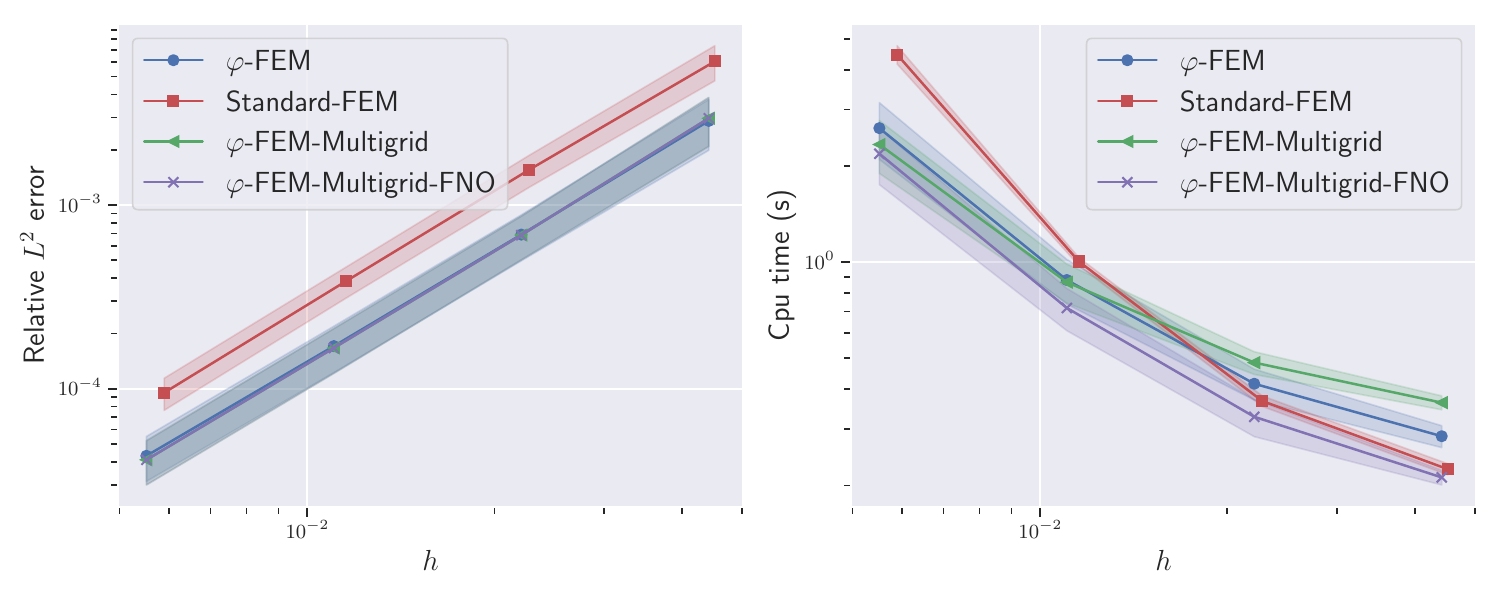}
      \caption{\textbf{Test case 3, $\mathbf{16 \times 16}$.} Left: relative $L^2$ errors with respect to the mesh size $h$. Right: CPU times (in seconds).}\label{fig:results_phi_fem_fno_multigrid_cas_test_1}
\end{figure}

The results in \cref{fig:results_phi_fem_fno_multigrid_cas_test_1} (left) illustrate that the three $\varphi$-FEM-based methods lead to comparable errors, both better than the Standard-FEM ones. However, in \cref{fig:results_phi_fem_fno_multigrid_cas_test_1} (right), we remark that these results are always obtained faster with $\varphi$-FEM-M-FNO. This gain is highlighted in \cref{table:phi_fem_m_fno_cas_test_1_16} where we have computed the average times and errors of each method for the five problems. Hence, we observe the gain on the coarse step, with comparable times for the steps. 

\begin{table}[!ht]
      \begin{tabular}{c|c|c|c|c|c}
            \multicolumn{1}{p{1.9cm}|}{\centering \vphantom{}                                                              \\ \textbf{Resolution} } &
            \multicolumn{1}{p{2.9cm}|}{\centering \vphantom{}                                                              \\ Method }&
            \multicolumn{1}{p{1.3cm}|}{\centering Time                                                                  \\ (coarse) }&
            \multicolumn{1}{p{1.3cm}|}{\centering Time                                                                  \\ (fine) }&
            \multicolumn{1}{p{1.3cm}|}{\centering Time                                                                  \\ (total) }&
            \multicolumn{1}{p{2.3cm}}{\centering $L^2$ relative \\ error}                                                        \\\hline\rule{0pt}{2.3ex}
            $\mathbf{32 \times 32}$   & Standard-FEM        & \diagbox{}{} & \diagbox{}{} & 0.23 & $6.06 \times 10^{-3}$ \\
                                      & $\varphi$-FEM       & \diagbox{}{} & \diagbox{}{} & 0.28 & $2.87 \times 10^{-3}$ \\
                                      & $\varphi$-FEM-M     & 0.16         & 0.20         & 0.36 & $2.96 \times 10^{-3}$ \\
                                      & $\varphi$-FEM-M-FNO & 0.004        & 0.21         & 0.21 & $2.96 \times 10^{-3}$ \\\hline\rule{0pt}{2.3ex}
            $\mathbf{64 \times 64}$   & Standard-FEM        & \diagbox{}{} & \diagbox{}{} & 0.36 & $1.56 \times 10^{-3}$ \\
                                      & $\varphi$-FEM       & \diagbox{}{} & \diagbox{}{} & 0.41 & $6.90 \times 10^{-4}$ \\
                                      & $\varphi$-FEM-M     & 0.16         & 0.32         & 0.48 & $6.85 \times 10^{-4}$ \\
                                      & $\varphi$-FEM-M-FNO & 0.004        & 0.32         & 0.32 & $6.85 \times 10^{-4}$ \\\hline\rule{0pt}{2.3ex}
            $\mathbf{128 \times 128}$ & Standard-FEM        & \diagbox{}{} & \diagbox{}{} & 1    & $3.85 \times 10^{-4}$ \\
                                      & $\varphi$-FEM       & \diagbox{}{} & \diagbox{}{} & 0.88 & $1.71 \times 10^{-4}$ \\
                                      & $\varphi$-FEM-M     & 0.16         & 0.70         & 0.86 & $1.67 \times 10^{-4}$ \\
                                      & $\varphi$-FEM-M-FNO & 0.004        & 0.70         & 0.70 & $1.67 \times 10^{-4}$ \\\hline\rule{0pt}{2.3ex}
            $\mathbf{256 \times 256}$ & Standard-FEM        & \diagbox{}{} & \diagbox{}{} & 4.46 & $9.51 \times 10^{-5}$ \\
                                      & $\varphi$-FEM       & \diagbox{}{} & \diagbox{}{} & 2.62 & $4.31 \times 10^{-5}$ \\
                                      & $\varphi$-FEM-M     & 0.16         & 2.17         & 2.33 & $4.1 \times 10^{-5}$  \\
                                      & $\varphi$-FEM-M-FNO & 0.004        & 2.17         & 2.18 & $4.1 \times 10^{-5}$
      \end{tabular}
      \caption{\textbf{Test case 3, $\mathbf{16 \times 16}$.} Results of the different methods. Both times and errors are the averages values on 5 data.}\label{table:phi_fem_m_fno_cas_test_1_16}
\end{table}

\paragraph*{$\mathbf{32 \times 32} $ grids}

\begin{figure}[!ht]
      \centering
      \includegraphics[width=0.8\textwidth]{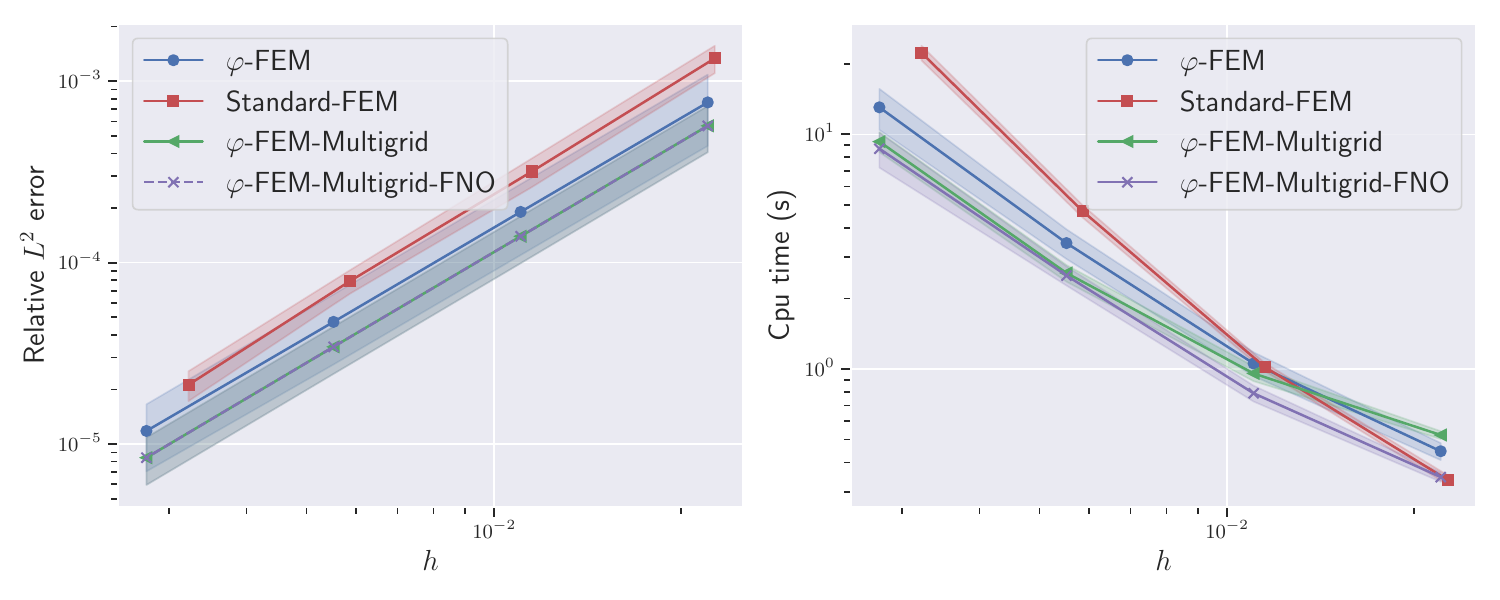}

      \caption{\textbf{Test case 3, $\mathbf{32 \times 32}$ data.} Left: relative $L^2$ errors with respect to the mesh size $h$. Right: CPU time with respect to $h$.}\label{fig:results_phi_fem_fno_multigrid_cas_test_1_32_32}
\end{figure}

Let us now consider an operator trained with data generated on finer grids, here $32 \times 32$, using the same amount of data and the same loss function. Once again, we consider five test data to compute the errors, for different fine resolutions, namely $64 \times 64$, $128\times 128$, $256 \times 256$, and $512 \times 512$. The results, given in \cref{fig:results_phi_fem_fno_multigrid_cas_test_1_32_32}, illustrate that using a better prediction to initialize the solver helps to reduce the computation time. Moreover, for a given tolerance of $10^{-9}$ as a stopping criterion, both multigrid approaches lead to better accuracy than the two other methods. In particular, we remark in \cref{table:phi_fem_m_fno_cas_test_1_32} that $\varphi$-FEM-M-FNO is always the fastest to compute, although the average computation time of the fine step is higher than when using $16 \times 16$ data, essentially due to the conversion between the Numpy matrix and the dolfinx vectors.

\begin{table}
      \begin{tabular}{c|c|c|c|c|c}
            \multicolumn{1}{p{1.9cm}|}{\centering \vphantom{}                                                              \\ \textbf{Resolution} } &
            \multicolumn{1}{p{2.9cm}|}{\centering \vphantom{}                                                              \\ Method }&
            \multicolumn{1}{p{1.3cm}|}{\centering Time                                                                  \\ (coarse) }&
            \multicolumn{1}{p{1.3cm}|}{\centering Time                                                                  \\ (fine) }&
            \multicolumn{1}{p{1.3cm}|}{\centering Time                                                                  \\ (total) }&
            \multicolumn{1}{p{2.3cm}}{\centering $L^2$ relative \\ error}                                                        \\\hline\rule{0pt}{2.3ex}
            $\mathbf{64 \times 64}$   & Standard-FEM        & \diagbox{}{} & \diagbox{}{} & 0.39  & $1.34 \times 10^{-3}$ \\
                                      & $\varphi$-FEM       & \diagbox{}{} & \diagbox{}{} & 0.50  & $7.64 \times 10^{-4}$ \\
                                      & $\varphi$-FEM-M     & 0.24         & 0.34         & 0.58  & $5.68 \times 10^{-4}$ \\
                                      & $\varphi$-FEM-M-FNO & 0.011        & 0.36         & 0.37  & $5.68 \times 10^{-4}$ \\\hline\rule{0pt}{3ex}
            $\mathbf{128 \times 128}$ & Standard-FEM        & \diagbox{}{} & \diagbox{}{} & 1.02  & $3.18 \times 10^{-4}$ \\
                                      & $\varphi$-FEM       & \diagbox{}{} & \diagbox{}{} & 1.09  & $1.90 \times 10^{-4}$ \\
                                      & $\varphi$-FEM-M     & 0.23         & 0.75         & 0.98  & $1.40 \times 10^{-4}$ \\
                                      & $\varphi$-FEM-M-FNO & 0.011        & 0.88         & 0.90  & $1.40 \times 10^{-4}$ \\
            \hline\rule{0pt}{3ex}
            $\mathbf{256 \times 256}$ & Standard-FEM        & \diagbox{}{} & \diagbox{}{} & 4.29  & $7.91 \times 10^{-5}$ \\
                                      & $\varphi$-FEM       & \diagbox{}{} & \diagbox{}{} & 3.46  & $4.71 \times 10^{-5}$ \\
                                      & $\varphi$-FEM-M     & 0.22         & 2.34         & 2.56  & $3.44 \times 10^{-5}$ \\
                                      & $\varphi$-FEM-M-FNO & 0.011        & 2.49         & 2.50  & $3.44 \times 10^{-5}$ \\\hline\rule{0pt}{3ex}
            $\mathbf{512 \times 512}$ & Standard-FEM        & \diagbox{}{} & \diagbox{}{} & 22.06 & $2.12 \times 10^{-5}$ \\
                                      & $\varphi$-FEM       & \diagbox{}{} & \diagbox{}{} & 13.16 & $1.18 \times 10^{-5}$ \\
                                      & $\varphi$-FEM-M     & 0.23         & 9.12         & 9.35  & $8.41 \times 10^{-6}$ \\
                                      & $\varphi$-FEM-M-FNO & 0.011        & 8.62         & 8.64  & $8.41 \times 10^{-6}$ \\
      \end{tabular}
      \caption{\textbf{Test case 3, $\mbf{32 \times 32}$ data.} Average computational time of the methods on 5 test data.}\label{table:phi_fem_m_fno_cas_test_1_32}
\end{table}

\subsubsection{Test case 4: the case of 3D complex geometries}

For this last test case, we consider only the three $\varphi$-FEM-based methods, not including the Standard-FEM results. The methods are used to solve the \cref{eq:main_poisson_equation_dirichlet_homo}, with non-homogeneous Dirichlet boundary conditions, on geometries defined by the level-set functions
\begin{equation*}
      \varphi(x,y,z) = (-1)^n \prod_j^n \bigg( -1 + \exp\left(- \frac{x_j^2}{2l_{x,j}^2} - \frac{y_j^2}{2l_{y,j}^2} -  \frac{z_j^2}{2l_{z,j}^2} \right)  \bigg) \,,
\end{equation*}
where
\[
      \begin{bmatrix}
            x_j \\
            y_j \\
            z_j
      \end{bmatrix}
      =
      R_z(\theta_z)\, R_y(\theta_y)\, R_x(\theta_x)\,
      \begin{bmatrix}
            x - \mu_x \\
            y - \mu_y \\
            z - \mu_z
      \end{bmatrix},
\]
with
\[
      R_x(\theta_x) =
      \begin{bmatrix}
            1 & 0            & 0             \\
            0 & \cos\theta_x & -\sin\theta_x \\
            0 & \sin\theta_x & \cos\theta_x
      \end{bmatrix}, 
      R_y(\theta_y) =
      \begin{bmatrix}
            \cos\theta_y  & 0 & \sin\theta_y \\
            0             & 1 & 0            \\
            -\sin\theta_y & 0 & \cos\theta_y
      \end{bmatrix}, 
      R_z(\theta_z) =
      \begin{bmatrix}
            \cos\theta_z & -\sin\theta_z & 0 \\
            \sin\theta_z & \cos\theta_z  & 0 \\
            0            & 0             & 1
      \end{bmatrix}\,,
\]
where all the parameters are randomly chosen, as in the previous test case. In addition, we consider 
\begin{equation*}
      f_{(A,\mu_0,\mu_1,\mu_2, \sigma_x, \sigma_y, \sigma_z)}(x,y,z) = A \exp\left(- \frac{(x-\mu_0)^2}{2\sigma_x^2} - \frac{(y-\mu_1)^2}{2\sigma_y^2} - \frac{(z-\mu_2)^2}{2\sigma_z^2}\right)\,,
\end{equation*}
and
\begin{equation*}
      g_{(\alpha, \beta)}(x,y) = \alpha \left((x-0.5)^2 - (y-0.5)^2\right) \cos\left( \beta z \pi \right)\,,
\end{equation*}
where all the parameters are randomly chosen following some uniform distributions. 

For the FNO-based approach, we generate a set of 250 data, divided into 200 training data and 50 validation data. These data are generated on Cartesian grids of resolution $20 \times 20 \times 20$, and we perform 200 training epochs, with batches of size 8. The loss function to minimize is the 3D extension of the $H^1$ approximation considered previously. 

\begin{figure}[!ht]
      \centering
      \includegraphics[width=0.8\textwidth]{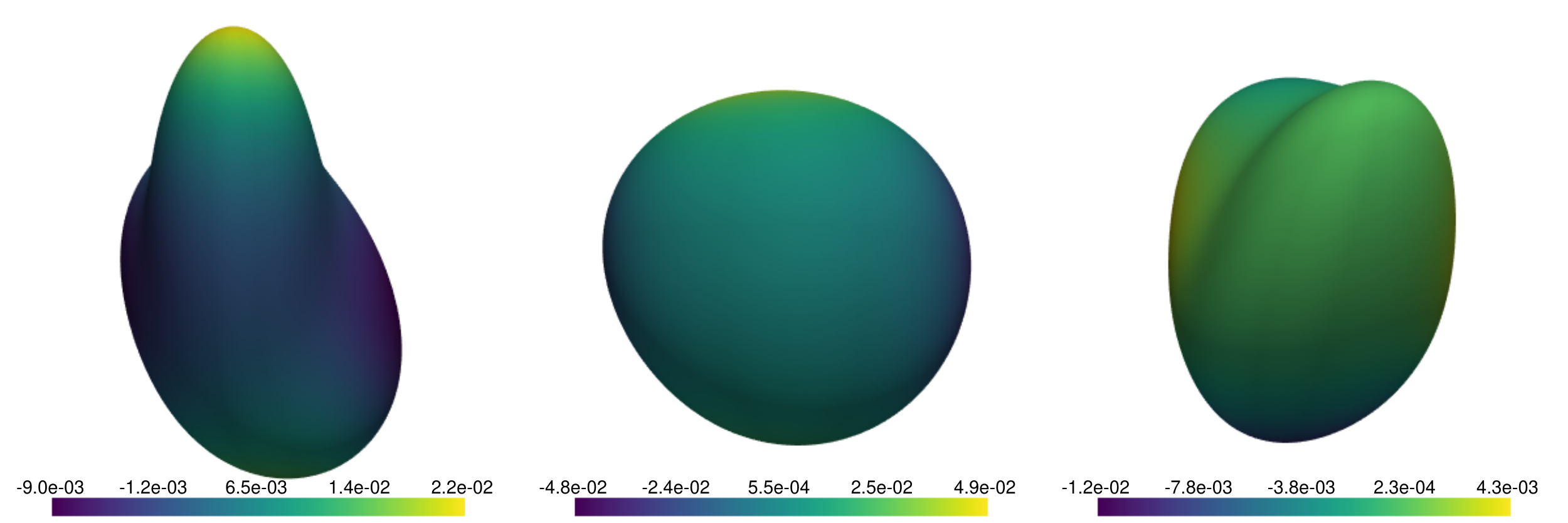}
      \caption{\textbf{Test case 4.} Representation of 3 reference solutions.}\label{fig:refrence_solutions_multigrid_fno_3D}
\end{figure}

Once the $\varphi$-FEM-FNO is trained, we consider a test sample of 6 new data to evaluate the performance of the methods. We represent 3 of the 6 reference solutions in \cref{fig:refrence_solutions_multigrid_fno_3D} to illustrate the variability of the considered geometries. The reference solutions are once again obtained by a fine conforming standard finite element method. 

Since we consider grids of size $20^3$ to generate data, the methods will be compared at 3 new resolutions: $40^3$, $80^3$, and $160^3$. For the three methods, we use once again the BiConjugate Gradient Stabilized with a stopping criterion set to $10^{-9}$. For $\varphi$-FEM-Multigrid, the coarse solution is obtained using a direct solver on grids of size $20^3$, $40^3$, and $80^3$. For the classical $\varphi$-FEM solver, we combine the iterative solver with an LU preconditioner.

\begin{figure}[!h]
      \centering
      \includegraphics[width=\twofigswidth]{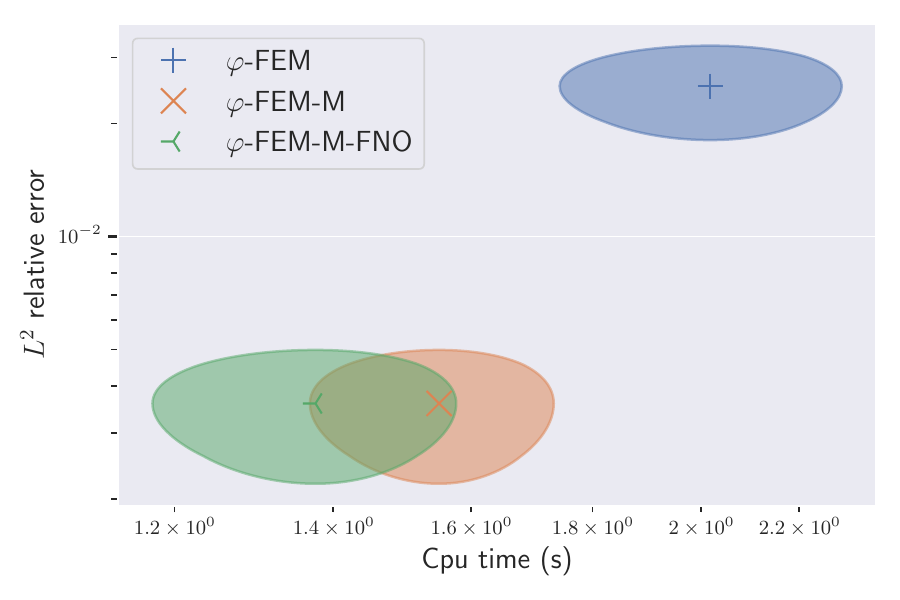}
      \quad
      \includegraphics[width=\twofigswidth]{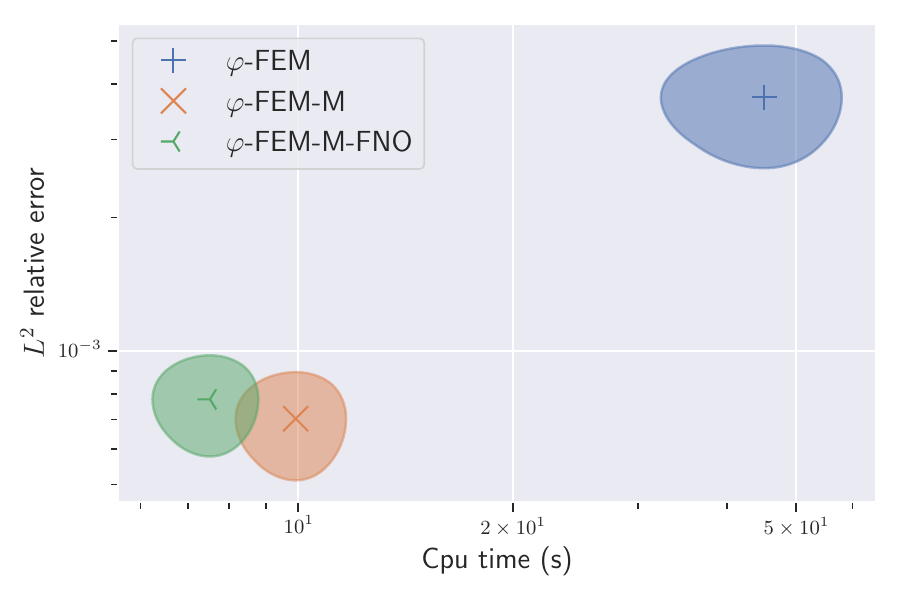}
      \\
      \includegraphics[width=\twofigswidth]{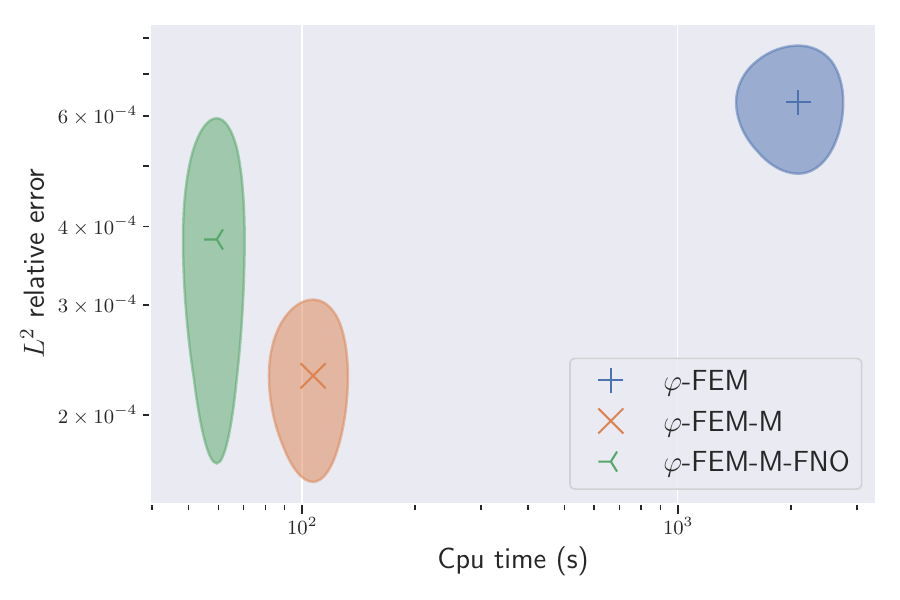}
      \caption{\textbf{Test case 4.} Results of the three methods. Top left: grids $40^3$. Top right: grids $80^3$. Bottom: grids $160^3$.}\label{fig:results_multigrid_fno_3D}
\end{figure}

To illustrate the performance of our multigrid approaches, we measure the relative $L^2$ errors and the computation times. The results given in \cref{fig:results_multigrid_fno_3D} present the average errors with respect to the average times (the colored areas corresponding to the standard deviation) and lead to the conclusion that $\varphi$-FEM-M-FNO is systematically much faster than the two approaches, even decreasing the error of the naive $\varphi$-FEM implementation.

\section{Conclusion}

The $\varphi$-FEM–M–FNO method represents a promising hybrid framework for the efficient resolution of partial differential equations posed on complex geometries. Indeed this new approach bridges classical scientific computing and machine learning techniques by combining:
\begin{itemize}
\item $\varphi$-FEM method, which accurately handles complex boundaries without requiring mesh conformity,
\item Multigrid methods, which ensure fast convergence and scalability across resolutions,
\item Fourier Neural Operators, which enable rapid inference once trained. This technique reduces computational cost while maintaining good accuracy especially in high-resolution.
\end{itemize}
Nonetheless, some challenges remain. The training phase of FNOs can be computationally expensive, and a theoretical analysis of convergence for the coupled numerical/learning framework is still an open problem.
This paper offers a first innovative foundation with strong potential for multi-scale and industrial applications where both accuracy and computational efficiency are critical.

\section*{Acknowledgments}
This work was supported by the Agence Nationale de la Recherche, Project PhiFEM, under grant ANR-22-CE46-0003-01. Vanessa Lleras thanks the project MSMPHI, under grant ANR-23-CE40-0013-01.

\bibliographystyle{abbrv}
\bibliography{biblio}

\end{document}